\documentclass[12pt]{amsart}
\usepackage{epsfig,a4wide}
\usepackage[english]{babel}
\usepackage[hypertex]{hyperref}
\usepackage{eepic}
\usepackage{color}

\newtheorem{theorem}{Theorem}[section]
\newtheorem{prop}[theorem]{Proposition}
\newtheorem{cor}[theorem]{Corollary}
\newtheorem{lemma}[theorem]{Lemma}

\newtheorem{rem}[theorem]{Remark}

\title{The oriented graph of multi-graftings in the Fuchsian case}
\author{Gabriel Calsamiglia}
\address{Instituto de Matem\'atica, Universidade Federal Fluminense,
Rua M\'ario Santos Braga s/n, 24020-140, Niter\'oi, Brazil} \email{gabriel@mat.uff.br}
\author{Bertrand Deroin}
\address{D\'epartement de  Math\'ematiques
d'Orsay, Universit\'e Paris 11, 91405 Orsay Cedex , France} \email{bertrand.deroin@math.u-psud.fr}
\author{Stefano Francaviglia}
\address{Dipartimento di
Matematica Universit\`a di Bologna, P.zza Porta S. Donato 5, 40126 Bologna, Italy } \email{stefano.francaviglia@unibo.it}

\thanks{B.D.'s research was partially supported by ANR-08-JCJC-0130-01,  ANR-09-BLAN-0116. G.C.'s was partially supported by MATH-AmSud and CNPq.}

\begin{document}
\begin{abstract}
We prove the connectedness and compute the diameter of the oriented graph of multi-graftings
associated to exotic
$\mathbb{CP}^1$-structures on a compact surface $S$ with a given holonomy representation of Fuchsian type.
\end{abstract}

\maketitle
\markright{The graph of multi-graftings in the Fuchsian case}

\section{Introduction}

Let $\Gamma_g$ be the fundamental group of a compact oriented surface $S$ of genus $g\geq 2$, and $\rho : \Gamma_g \rightarrow \mathrm{PSL} (2,\mathbb R)$ be a \textit{Fuchsian} representation, namely a faithful and discrete one. A marked surface of genus $g$ is the data of a simply connected cover $\widetilde{S} $ of $S$ together with a free discontinuous action of $\Gamma_g$. A $\mathbb {CP}^1$-structure (sometimes referred to as a projective structure) with holonomy $\rho$ on the marked surface is a local diffeomorphism $D : \widetilde{S} \rightarrow \mathbb {CP}^1 $ called developing map which is $\rho$-equivariant. We denote by $P(\rho)$ the set of equivalence classes of marked $\mathbb{ C P}^1$-structures on a surface of genus $g$ with holonomy $\rho$, where two projective structures $(\widetilde{S_i}, D_i)$, $i=1,2$ are equivalent if there exists a $\Gamma_g$-equivariant diffeomorphism $\Phi : \widetilde{S_1} \rightarrow \widetilde{S_2}$ such that $D_1 = D_2 \circ \Phi$.\footnote{This definition of projective structure coincides with the classical one because there is no ambiguity in the choice of developing map when the holonomy representation is non-elementary, see ~\cite[Lemma 12.10]{CDF}.}

This article deals with the study of a surgery operation called
\textit{grafting} that produces, given an element in $P (\rho)$, new elements in the same set. Grafting consists in cutting a surface
equipped with a $\mathbb{CP}^1$-structure along a particular type of simple closed curve called
\textit{graftable curve}, and gluing a Hopf annulus, namely the quotient of a simply connected
domain of the Riemann sphere invariant by the (loxodromic) holonomy of the graftable
curve. This operation produces a new element of $P(\rho)$.

 Grafting was used by Hejhal \cite[Theorem 4]{Hej} and Thurston (unpublished)  to produce examples of projective structures with holonomy $\rho$ that are different from the uniformizing structure $\sigma_u= \rho(\Gamma_g) \backslash \mathbb H^2$. Such structures are called \textit{exotic}. The importance of grafting comes from the fact that it allows to define coordinates on $\mathcal{P}(\rho)$ when $\rho$ is a Fuchsian representation: Goldman proved that any $\mathbb{CP}^1$-structure with holonomy $\rho$ is obtained from the uniformizing one by grafting a collection of disjoint graftable simple closed curves (see~\cite{Goldman1}). Such an operation will be called a \textit{multi-grafting}.

The goal of this note is to improve Goldman's result in the following way.

\begin{theorem} \label{t:positive connectivity}  Let $\sigma_1$ and $\sigma_2$ be two exotic
projective structures sharing the same Fuchsian holonomy. Then
$\sigma_2$ can be obtained from
$\sigma_1$ by a sequence of two multi-graftings. \end{theorem}

A consequence of this result is that there exist positive
cycles of graftings, namely finite sequences of marked $\mathbb{CP}^1$-structures $\sigma_0,\ldots,\sigma_r=\sigma_0$ such that for each $i=1,\ldots, r$,
$\sigma_i$ is a grafting of $\sigma_{i-1}$. The integer $r$ is then called the period of the cycle. Observe that an immediate corollary of the theorem is that
any couple of exotic $\mathbb{CP}^1$-structures are contained in such
a positive cycle of period bounded by~$4$.
We will see (Corollary~\ref{l2c})
that indeed there are such cycles of period~$2$.

Let  $MG (\rho)$ be the oriented graph whose vertices are elements of $\mathcal{P}(\rho)$ and two vertices $\sigma_1,\sigma_2$ are joined by an oriented edge from $\sigma_1$ to $\sigma_2$ if $\sigma_2$ is obtained from $\sigma_1$ by a multi-grafting.  Theorem \ref{t:positive connectivity} can be restated by saying that the oriented graph of multi-graftings $MG(\rho)\setminus\sigma_u$ is a connected graph of radius $2$. As a consequence we also get that the fundamental group of $MG(\rho)$ is not finitely generated.

To prove the results we will use some surgery operations on multi-curves introduced by
Luo~\cite{L} and later developed by Ito~\cite{ITO}. Our results and methods are closely related
to Thompson's, see~\cite{Thompson}, but he considers the case of Schottky representations
instead of Fuchsian ones. We observe that our argument extends {\em stricto sensu} to the case of quasi-Fuchsian representations.

\section{Graftable curves}
In this section we introduce the action of grafting on $\mathcal{P}(\rho)$ and define the graph of multi-graftings.

\subsection{Definition}
Recall that a \textit{multi-curve} on a surface $S$ is a finite disjoint union of simple closed curves none of which is homotopically
trivial. Let $\sigma$ be a marked projective structure on a compact orientable surface $S$. A multi-curve is said to be \textit{graftable} (in $\sigma$) if all of its
components have loxodromic holonomy and the developing map is injective when restricted to a lift of any of those components in $\widetilde{S}$. The condition is
independent of the choice of representative in the class $[\sigma]\in\mathcal{P}(\rho)$.

\subsection{Grafting along graftable curves}
If $\alpha = \{\alpha_i\}_{i\in I}$ is a graftable multi-curve, one can produce  another marked projective structure,
called the grafting along $\alpha$, and denoted $\text{Gr}(\sigma,\alpha)$. We recall the construction here. We cut the surface $\widetilde{S}$ along the lifts
$\widetilde{\alpha_i}$'s of the curves $\alpha_i$'s, and glue to each of them a copy of $\mathbb{CP}^1\setminus \overline{D(\widetilde{\alpha_i}})$ using the
developing map for the gluing. We then obtain a new surface denoted by $\widetilde{S'}$, together with a new map $D' : \widetilde{S'} \rightarrow \mathbb{CP}^1$
which is defined by $D$ on $\widetilde{S}\setminus \pi^{-1} (\cup_i \alpha_i)$ and by the identity on the spherical domains $\mathbb{CP}^1 \setminus \overline{D
  (\widetilde{\alpha_i})}$. The $\Gamma_g$-action on $\widetilde{S}$
induces a $\Gamma_g$-action on $\widetilde{S'}$ which is free and discontinuous, and the map $D'$ is obviously $\rho$-equivariant. Hence, this defines a new
marked projective structure $\text{Gr} (\sigma , \alpha)$ with holonomy $\rho$: the grafting of $\sigma$ over the graftable multi-curve $\alpha$.

As  $\alpha_i$ has loxodromic holonomy, it acts freely and properly discontiuosly on $\mathbb{CP}^1\setminus\overline{D(\widetilde{\alpha_i})}$, and its quotient
is a cylinder equipped with a projective structure. Therefore, the grafting can be viewed as a cut-and-paste procedure directly in $S$, which cuts $S$ along each
$\alpha_i$ and glues back the cylinder
$\langle \alpha_i\rangle\backslash(\mathbb{CP}^1\setminus\overline{D(\widetilde{\alpha_i})})$.

\subsection{Isotopy class of graftable curves}\label{ss:isotopy}
It is an easy fact to verify that if $\alpha$ and $\alpha'$ belong to the same connected
component of the set of \emph{graftable} multi-curves (for the compact open topology), then the resulting projective structures $\text{Gr} (\sigma , \alpha)$ and
$\text{Gr} (\sigma , \alpha')$ are equivalent. However, we will see that it can happen that $\alpha$ and $\alpha'$ are two graftable multi-curves that are isotopic as
multi-curves by an isotopy that leaves the space of graftable
multi-curves, and such that their corresponding graftings  are not
equivalent (see Remark~\ref{RBRBL}).

\subsection{The graph of multi-graftings}\label{ss:graph}
Let $\rho$ be a representation from $\Gamma_g$ to $\mathrm{PSL}(2,\mathbb C)$. Let us define the graph
of multi-graftings $MG(\rho)$ in the following way. The vertices are the elements of $\mathcal{P}(\rho)$ and two of them  $(S_1,\sigma_1)$ and $(S_2,\sigma_2 )$ are the connected by a positive segment from $\sigma_1$ to $\sigma_2$ if there exists a graftable multi-curves $\alpha$ in $S_1$ such that $\mathrm{Gr} (\sigma_1, \alpha) = \sigma_2$.

\section{Fuchsian case: construction of graftable curves}

Recall that a representation $\rho : \Gamma_g \rightarrow \mathrm{PSL} (2,\mathbb R)$ is Fuchsian if it is discrete and faithful. In the sequel $\rho$ will
always be assumed to be Fuchsian.

\subsection{Goldman's parametrization of
  $MG(\rho)$}\label{ss:Goldman's parametrization}
We will denote by $\sigma_u$ the uniformizing structure on the surface $S_u := \rho(\Gamma_g) \backslash \mathbb H^2$, which is obtained by taking the quotient
of $\mathbb H^2$ by the $\rho$-action of $\Gamma_g$ on $\mathbb H^2$. For this structure, the developing map is just the identity when identifying the universal
cover of $S_u$ with $\mathbb H^2$, and in particular is injective. Hence, any simple closed curve on $S_u$ is a graftable curve.  Hence in this case the space of graftable multi-curves and the space of multi-curves are the same. By the discussion in~\S\ref{ss:isotopy}  the
 grafting $\text{Gr} (\sigma_u, \alpha)$ depends only on the isotopy class of $\alpha$ as a \emph{multi-curve}.

Goldman proved in \cite{Goldman1} that every marked projective structure $\sigma$ with holonomy $\rho$ is obtained by grafting the structure $\sigma_u$ along a multi-curve $\alpha
= \{ \alpha_i \}_i$. Moreover, this family is unique, and can be reconstructed from $\sigma$ in the following way. For a Fuchsian  projective structure $\sigma$,
denote by $S^{\mathbb R}$ (resp. $S^{\pm}$) the quotient of $D^{-1} (\mathbb{RP}^1)$ (resp. $D^{-1} (\mathbb H^{\pm})$) by the covering group $\Gamma_g$. Since
$\rho$ is Fuchsian, it preserves the decomposition $\mathbb {CP}^1 = \mathbb H^+ \cup \mathbb{RP}^1 \cup \mathbb H^-$, and thus $S^{\mathbb R}$ is an analytic
real submanifold of $S$ separating $S$ in domains which are either positive or negative according they belong to $S^+$ or $S^-$. Goldman proved that the
components of $S^-$ are necessarily annuli. The set of annuli is homotopic to a unique multi-loop $\alpha$ satisfying $\sigma =
\text{Gr}(\sigma_u, \alpha)$.  To abridge notations we define $\text{Gr}_\alpha := \text{Gr}(\sigma_u,\alpha)$.

\subsection{Homotopically transverse multi-curves} Let $\alpha = \{ \alpha_i\}_{i\in I}$ and $\beta = \{ \beta_j \}_{j\in J}$ be two multi-curves.  They are
homotopically transverse if the following conditions hold: \begin{itemize} \item for each $i\in I$ and $j\in J$, the curves $\alpha_i$ and $\beta_j$ are not
homotopic \item  they are transverse in the usual sense and
\item The complement of $(\cup \alpha_i)\cup(\cup\beta_j)$ in $S$ has no bi-gon component. 
%any component of $\beta_j\cap(S\setminus\cup_i\alpha_i)$ sitting in a component $C$ of
%$S\setminus\cup_i\alpha_i$ is not homotopically equivalent
%    (fixing endpoints) in $\overline{C}$ to a segment in $\partial C$ and the same remains true if we exchange the roles of the multi-curves $\alpha$ and
%    $\beta$.
\end{itemize}

\subsection{Construction of graftable multi-curves}\label{ss:construction}

Given a multi-curve $\alpha=\{\alpha_i\}_{i\in I}$, a set of turning directions for $\alpha$ is an assignment to each curve $\alpha_i$ of a turning direction
$T_i\in\{R,L\}$ (``Right'' or ``Left'') in such a way that any two parallel curves have the same turning direction.

In this paragraph we provide a construction that, given two homotopically transverse multi-curves $\alpha=\{\alpha_i\}_{i\in I}$  and  $\beta=\{\beta_j\}_{j\in
J}$, and a set $T=\{T_i\}$ of turning directions for $\alpha$, produces a multi-curve $\beta_{T}$ which is graftable in $\text{Gr}_\alpha$ and isotopic to
$\beta$.

We begin by assuming that there are no parallel curves in the families $\alpha$ and $\beta$. In this case we can assume that the components of $\alpha$ and
$\beta$ are simple closed geodesics  in the uniformizing structure $\sigma_u$.

Recall that $\text{Gr}_{\alpha}$ is obtained by gluing $S_u\setminus\alpha$ with some grafting annuli. We will explain the construction  of $\beta_T$ in each piece of this decomposition separately, beginning with the intersection of $\beta_T$ with $S_u\setminus\alpha$, and then construct the intersection of $\beta_{T}$ with the grafting
annuli glued to $S_u\setminus\alpha$ to obtain $\text{Gr}_\alpha$.

The boundary of $S_u \setminus \alpha$ consists of two copies $\alpha_i'$ and $\alpha_i''$ of each curve $\alpha_i$, and for each component $C$ of $S_u \setminus
\alpha$, its boundary is a union of such components. We fix a small positive number $\varepsilon$, and for each $p\in\alpha_i\cap\beta\in\partial C$, we consider
the point $p_T\in\partial C$ lying at distance $\varepsilon$ from $p$ to the side of $p$ indicated by $T_i$  with respect to the orientation induced on
$\alpha_i$ by $C$. If we do this for all components of $S_u\setminus\alpha$, we get for each point $p\in\alpha_i\cap\beta$ a couple of distinct points
$p'\in \alpha_i '$ and $p'' \in \alpha_i ''$ lying at distance $\varepsilon$ from $p$ (as seen as a point in $\alpha_i$ or $\alpha_i ''$ under the natural identifications $\alpha_i \simeq \alpha_i ' \simeq \alpha_i''$).

Now, $\beta\cap C$ is a union of geodesic segments $[p,q]$ joining points of $\partial C$. We define $\beta_{T}$ in $S_u \setminus \alpha\subset S$ to be the
union of the segments $[p_T,q_T]$ with $p_T$ and $q_T$ constructed as above. Observe that if we move the points $p,q$ a little bit, then the segments $[p_T,q_T]$
are disjoint in the component $C$, but also in the whole surface $S$.

Then, one has to define the curve $\beta_{T}$ in the grafting annuli in a graftable way. The continuation should start from the point $p'$ above and end at
$p''$. (In Figure~\ref{f0} we depicted the case $T_i=L$.)

\setlength{\unitlength}{1ex} 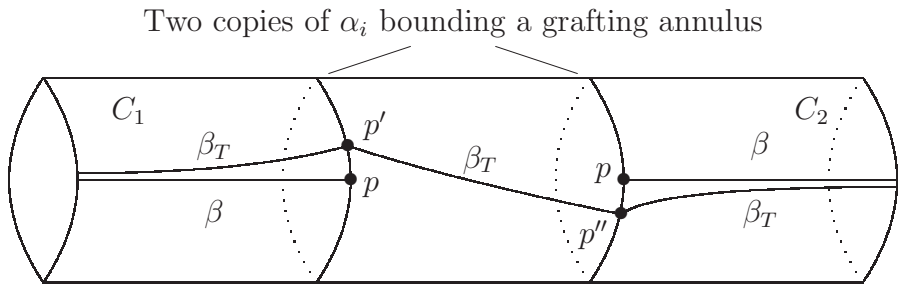
\begin{figure}[htbp]
  \centering
  \begin{picture}(60,20)
\multiput(0,0)(0,15){2}{\line(1,0){60}} \multiput(0,0)(20,0){4}{\qbezier(0,0)(5,7.5)(0,15)} \multiput(20,0)(20,0){3}{\bezier{20}(0,0)(-5,7.5)(0,15)}
\qbezier(0,0)(-5,7.5)(0,15)

\put(22.5,7.5){\makebox(0,0){$\bullet$}} \put(22.5,7.5){\makebox(3,-1){$p$}} \put(42.5,7.5){\makebox(0,0){$\bullet$}} \put(42.5,7.5){\makebox(-3,1){$p$}}
\put(22.3,10){\makebox(0,0){$\bullet$}} \put(22.3,10){\makebox(4,3){$p'$}} \put(42.3,5){\makebox(0,0){$\bullet$}} \put(42.3,5){\makebox(-4,-3){$p''$}}

\multiput(2.5,7.5)(40,0){2}{\line(1,0){20}} \multiput(12.5,5)(40,5){2}{\makebox(0,0){$\beta$}} \qbezier(2.5,8)(15,8)(22.3,10) \qbezier(22.3,10)(30,7.5)(42.3,5)
\qbezier(42.3,5)(45,7)(62.5,7) \multiput(12.5,10)(40,-5){2}{\makebox(0,0){$\beta_T$}} \put(32,9){\makebox(0,0){$\beta_T$}}

\put(20.9,15.3){\line(3,1){6}} \put(39.1,15.3){\line(-3,1){6}} \put(30,19){\makebox(0,0){Two copies of $\alpha_i$ bounding a grafting annulus}}

\put(5,12){$C_1$} \put(55,12){$C_2$}
  \end{picture}
  \caption{The curve $\beta_T$ in the surface $S_u$. Here $\alpha_i$
    appears in the boundary of two components $C_1$ and $C_2$. In the
    picture, we used $T_i=L$.}
  \label{f0}
\end{figure}

To be sure that $\beta_T$ is graftable and in the isotopy class of $\beta$, we need some care. First, we suppose that $\beta$ intersects $\alpha_i$ once. Figure
\ref{fig:imagebeta} provides a sketch of the construction in the universal cover (we used the convention that $\mathbb H^2=\mathbb H^+$ is the upper
half-plane.)

\setlength{\unitlength}{1ex}
\begin{figure}[htbp]
  \centering
  \begin{picture}(80,35)
\put(0,-5){
    \multiput(0,20)(44,0){2}{\line(1,0){36}}
    \multiput(18.5,20)(43,0){2}{\line(0,1){20}}
    \multiput(16,38)(46,0){2}{$\widetilde\beta$}
    \multiput(18,20)(44,0){2}{\arc{24}{3.14}{0}}
    \multiput(-1,21)(44,0){2}{\footnotesize $\mathbb H^+$}
    \multiput(-1,17)(44,0){2}{\footnotesize $\mathbb H^-$}
{\color{red}
            \put(14,20){\arc{24}{0}{4.89}}
            \put(10,20){\arc{32}{-.805}{0}}
            \put(66,20){\arc{24}{-1.75}{3.14}}
            \put(70,20){\arc{32}{3.14}{3.94}}

            \qbezier(21,31.5)(19,35)(19,40)
            \qbezier(16,32)(18,30)(18,25)
            \qbezier(64,32)(62,30)(62,25)
            \qbezier(59,31.5)(61,35)(61,40)

            \multiput(21,37)(33,0){2}{$D(\widetilde\beta_T)$}
            \multiput(22,25)(52,5){2}{$\widetilde{\eta}'$}
            \multiput(4,29)(52,-4){2}{$\widetilde{\eta''}$}
            \multiput(13,10)(53,0){2}{$\widetilde\eta$}
          }

    \put(8,24){\footnotesize$\widetilde\alpha_i$}
    \put(49,25){\footnotesize$\widetilde\alpha_i$}

    \qbezier(9.6,31)(10,30)(10.5,29.2)
    \qbezier(70.4,31)(70,30)(69.5,29.2)
    \qbezier(25.5,29.2)(25,28.5)(23.8,28)
    \qbezier(54.5,29.2)(55,28.5)(56.2,28)

    \put(11,31){\makebox(0,-1){\footnotesize$\varepsilon$}}
    \put(24,29){\makebox(0,1){\footnotesize$\varepsilon$}}
    \put(56,29){\makebox(1,1){\footnotesize$\varepsilon$}}
    \put(69,30){\makebox(.5,1.2){\footnotesize$\varepsilon$}}

    \put(16,32){\makebox(0,-.5){\footnotesize$\bullet$}}
    \put(16,32){\makebox(0,3){\footnotesize$\widetilde{p''}$}}
    \put(59,31){\makebox(0,1){\footnotesize$\bullet$}}
    \put(59,31){\makebox(-3,4){\footnotesize$\widetilde{p'}$}}
    \put(21,31){\makebox(0,1){\footnotesize$\bullet$}}
    \put(21,31){\makebox(3,5){\footnotesize$\widetilde{p'}$}}
    \put(64,32){\makebox(0,-.5){\footnotesize$\bullet$}}
    \put(64,32){\makebox(0,3){\footnotesize$\widetilde{p''}$}}

    \put(2.2,20){\makebox(0,0){\footnotesize$\bullet$}}
    \put(4,22){\makebox(0,0){\footnotesize$\widetilde{\xi''}$}}

    \put(54.2,20){\makebox(0,0){\footnotesize$\bullet$}}
    \put(56,22){\makebox(0,0){\footnotesize$\widetilde{\xi''}$}}

    \put(26,20){\makebox(0,0){\footnotesize$\bullet$}}
    \put(27,18){\makebox(0,0){\footnotesize$\widetilde{\xi'}$}}

    \put(78,20){\makebox(0,0){\footnotesize$\bullet$}}
    \put(79,18){\makebox(0,0){\footnotesize$\widetilde{\xi'}$}}
}
    \put(10,0){\footnotesize Case $T_i=L$}
    \put(60,0){\footnotesize Case $T_i=R$}
  \end{picture}

  \caption{The portion of $\beta_T$ in the universal cover of the
    grafting annulus}
  \label{fig:imagebeta}
\end{figure}
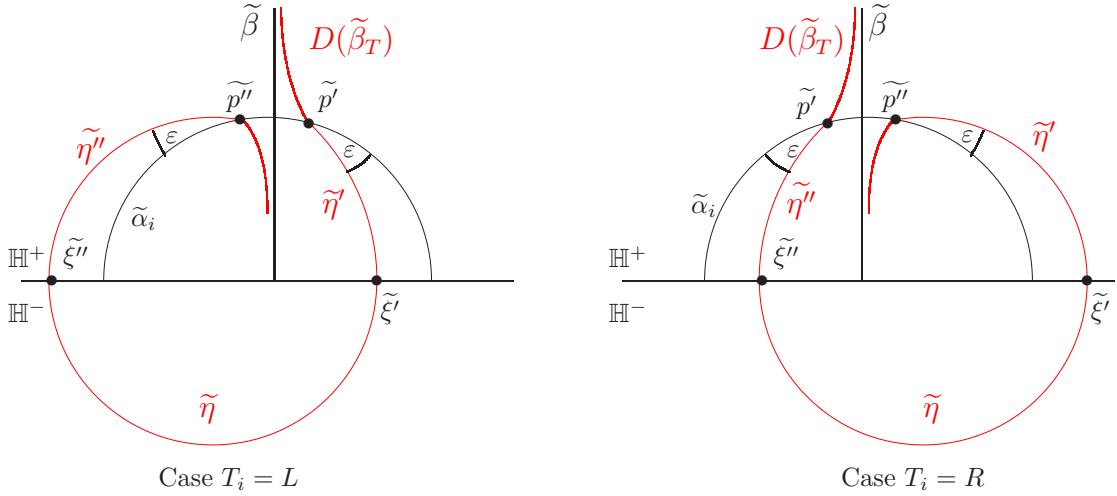

When the path $\beta_{T}$ enters in the grafting, it means that any lift $\widetilde{\beta_{T}}$ enters in the subset $\mathbb {CP}^1\setminus
\overline{\widetilde{\alpha_i}}$ that we have glued to $\widetilde{S_u}$ to obtain $\widetilde{\text{Gr}_\alpha}$. It enters at the point $\widetilde{p}'$ and needs to get
out at the point $\widetilde{p}''$ by a path in $\mathbb {CP}^1\setminus \overline{\widetilde{\alpha_i}}$. For this it has to turn around the segment
$\overline{\widetilde{\alpha_i}}$ in the sphere. Since we want a graftable curve we need to avoid creating self-intersection points of the developed image of
$\widetilde\beta_T$. An example of such a curve can be constructed as follows. Consider two semi-infinite geodesics $\widetilde{\eta'}$ and $\widetilde{\eta''}$
starting from $\widetilde{p'}$ and $\widetilde{p''}$ and forming and angle $\varepsilon$ with $\widetilde\alpha_i$ as in Figure~\ref{fig:imagebeta}. Such
geodesics meet the real line (i.e. the boundary of $\mathbb H^+$) at two points $\widetilde{\xi'},\widetilde{\xi''}$.

When $\widetilde\beta_T$ meets $\widetilde\alpha_i$ at the point $\widetilde{p'}$, we continue it by $\widetilde{\eta'}$, then in $\mathbb H^-$ by the geodesic
$\widetilde\eta$ between $\widetilde{\xi'}$ and $\widetilde{\xi''}$, and finally with $\widetilde{\eta''}$. (See Figure~\ref{fig:imagebeta}.)

The path $\widetilde{\eta'}\ast\widetilde{\eta}\ast\widetilde{\eta''}$ takes values in the set $\mathbb{CP}^1\setminus\overline{\widetilde\alpha_i}$. Such path remains
embedded when quotienting $\mathbb{CP}^1\setminus\overline{\widetilde\alpha_i}$ by the action of $\alpha_i$ and provides the path $\beta_T$ in the grafting
annulus. Moreover, since $\mathbb{CP}^1\setminus\overline{\widetilde\alpha_i}$ is a disc, any two paths joining two points in the boundary are homotopic. This
shows that $\beta_T$ is indeed isotopic to $\beta$. (See also Figure~\ref{f4}.)

Let us do the construction when
 $\beta$ intersects $\alpha_i$ in more than one point. What we need to
 describe  is the part of $\beta_T$ in the grafting annulus. Again, we
 work in the universal cover. In Figure~\ref{fig:morethanone} we sketched
 the case of two points of intersection.

\setlength{\unitlength}{1ex}
\begin{figure}[htbp]
  \centering
  \begin{picture}(80,35)
\put(0,-5){
    \multiput(0,20)(44,0){2}{\line(1,0){36}}
    \multiput(18,20)(44,0){2}{\arc{24}{3.14}{0}}
    \multiput(-1,21)(44,0){2}{\footnotesize $\mathbb H^+$}
    \multiput(-1,17)(44,0){2}{\footnotesize $\mathbb H^-$}

{\color{red}
% HERE a dashed arc equivalent to
%            \put(15.85,20){\arc{22}{0}{4.34}}
%
            \put(15.85,20){\arc{22}{0}{.1}}
            \put(15.85,20){\arc{22}{.2}{.3}}
            \put(15.85,20){\arc{22}{.4}{.5}}
            \put(15.85,20){\arc{22}{.6}{.7}}
            \put(15.85,20){\arc{22}{.8}{.9}}
            \put(15.85,20){\arc{22}{1}{1.1}}
            \put(15.85,20){\arc{22}{1.2}{1.3}}
            \put(15.85,20){\arc{22}{1.4}{1.5}}
            \put(15.85,20){\arc{22}{1.6}{1.7}}
            \put(15.85,20){\arc{22}{1.8}{1.9}}
            \put(15.85,20){\arc{22}{2}{2.1}}
            \put(15.85,20){\arc{22}{2.2}{2.3}}
            \put(15.85,20){\arc{22}{2.4}{2.5}}
            \put(15.85,20){\arc{22}{2.6}{2.7}}
            \put(15.85,20){\arc{22}{2.8}{2.9}}
            \put(15.85,20){\arc{22}{3}{3.1}}
            \put(15.85,20){\arc{22}{3.2}{3.3}}
            \put(15.85,20){\arc{22}{3.4}{3.5}}
            \put(15.85,20){\arc{22}{3.6}{3.7}}
            \put(15.85,20){\arc{22}{3.8}{3.9}}
            \put(15.85,20){\arc{22}{4}{4.1}}
            \put(15.85,20){\arc{22}{4.2}{4.34}}

% HERE a dashed arc equivalent to
%          \put(14.85,20){\arc{24}{-1.5}{0}}
%
            \put(14.85,20){\arc{24}{-1.5}{-1.4}}
            \put(14.85,20){\arc{24}{-1.3}{-1.2}}
            \put(14.85,20){\arc{24}{-1.1}{-1}}
            \put(14.85,20){\arc{24}{-.9}{-.8}}
            \put(14.85,20){\arc{24}{-.7}{-.6}}
            \put(14.85,20){\arc{24}{-.5}{-.4}}
            \put(14.85,20){\arc{24}{-.3}{-.2}}
            \put(14.85,20){\arc{24}{-.1}{0}}
% HERE a dashed arc equivalent to
%            \put(63.85,20){\arc{26}{0}{4.05}}
%
\put(48,0){ \put(15.85,20){\arc{26}{0}{.1}}
            \put(15.85,20){\arc{26}{.2}{.3}}
            \put(15.85,20){\arc{26}{.4}{.5}}
            \put(15.85,20){\arc{26}{.6}{.7}}
            \put(15.85,20){\arc{26}{.8}{.9}}
            \put(15.85,20){\arc{26}{1}{1.1}}
            \put(15.85,20){\arc{26}{1.2}{1.3}}
            \put(15.85,20){\arc{26}{1.4}{1.5}}
            \put(15.85,20){\arc{26}{1.6}{1.7}}
            \put(15.85,20){\arc{26}{1.8}{1.9}}
            \put(15.85,20){\arc{26}{2}{2.1}}
            \put(15.85,20){\arc{26}{2.2}{2.3}}
            \put(15.85,20){\arc{26}{2.4}{2.5}}
            \put(15.85,20){\arc{26}{2.6}{2.7}}
            \put(15.85,20){\arc{26}{2.8}{2.9}}
            \put(15.85,20){\arc{26}{3}{3.1}}
            \put(15.85,20){\arc{26}{3.2}{3.3}}
            \put(15.85,20){\arc{26}{3.4}{3.5}}
            \put(15.85,20){\arc{26}{3.6}{3.7}}
            \put(15.85,20){\arc{26}{3.8}{3.9}}
            \put(15.85,20){\arc{26}{4}{4.05}}
}
% HERE a dashed arc equivalent to
%            \put(64.25,20){\arc{25.2}{-1.9}{0}}
%
            \put(64.25,20){\arc{25.2}{-1.9}{-1.8}}
            \put(64.25,20){\arc{25.2}{-1.7}{-1.6}}
            \put(64.25,20){\arc{25.2}{-1.5}{-1.4}}
            \put(64.25,20){\arc{25.2}{-1.3}{-1.2}}
            \put(64.25,20){\arc{25.2}{-1.1}{-1}}
            \put(64.25,20){\arc{25.2}{-.9}{-.8}}
            \put(64.25,20){\arc{25.2}{-.7}{-.6}}
            \put(64.25,20){\arc{25.2}{-.5}{-.4}}
            \put(64.25,20){\arc{25.2}{-.3}{-.2}}
            \put(64.25,20){\arc{25.2}{-.1}{0}}

            \bezier{30}(12,30)(14,30)(15,25)
            \bezier{30}(16,32)(14,34)(13,37)

\color{blue}
% HERE a dashed arc equivalent to
%            \put(15.85,20){\arc{25.5}{0}{5.1}}
%
            \put(15.85,20){\arc{25.5}{0}{0.2}}
            \put(15.85,20){\arc{25.5}{.4}{.6}}
            \put(15.85,20){\arc{25.5}{.8}{1}}
            \put(15.85,20){\arc{25.5}{1.2}{1.4}}
            \put(15.85,20){\arc{25.5}{1.6}{1.8}}
            \put(15.85,20){\arc{25.5}{2}{2.2}}
            \put(15.85,20){\arc{25.5}{2.4}{2.6}}
            \put(15.85,20){\arc{25.5}{2.8}{3}}
            \put(15.85,20){\arc{25.5}{3.2}{3.4}}
            \put(15.85,20){\arc{25.5}{3.6}{3.8}}
            \put(15.85,20){\arc{25.5}{4}{4.2}}
            \put(15.85,20){\arc{25.5}{4.4}{4.6}}
            \put(15.85,20){\arc{25.5}{4.8}{5.1}}

% HERE a dashed arc equivalent to
%            \put(14,20){\arc{29.2}{-.74}{0}}
%
            \put(14,20){\arc{29.2}{-.74}{-.54}}
            \put(14,20){\arc{29.2}{-.34}{-.14}}

% HERE a dashed arc equivalent to
%            \put(64,20){\arc{22}{-1.1}{3.14}}
%
            \put(64,20){\arc{22}{-1.1}{-.9}}
            \put(64,20){\arc{22}{-.7}{-.5}}
            \put(64,20){\arc{22}{-.3}{-.1}}
            \put(64,20){\arc{22}{.1}{.3}}
            \put(64,20){\arc{22}{.5}{.7}}
            \put(64,20){\arc{22}{.9}{1.1}}
            \put(64,20){\arc{22}{1.3}{1.5}}
            \put(64,20){\arc{22}{1.7}{1.9}}
            \put(64,20){\arc{22}{2.1}{2.3}}
            \put(64,20){\arc{22}{2.5}{2.7}}
            \put(64,20){\arc{22}{2.9}{3.14}}

% HERE a dashed arc equivalent to
%            \put(64.7,20){\arc{23.4}{3.14}{4.7}}
%
            \put(64.7,20){\arc{23.4}{3.14}{3.34}}
            \put(64.7,20){\arc{23.4}{3.54}{3.74}}
            \put(64.7,20){\arc{23.4}{3.94}{4.14}}
            \put(64.7,20){\arc{23.4}{4.34}{4.54}}

            \bezier{20}(21,31.6)(23,31)(24,25)
            \bezier{20}(24.8,29.7)(23,31)(24,38)
          }

    \put(8,24){\footnotesize$\widetilde\alpha_i$}
    \put(49,25){\footnotesize$\widetilde\alpha_i$}

    \put(16,32){\makebox(0,-.5){\footnotesize$\bullet$}}
    \put(11,31){\makebox(0,3){\footnotesize$\widetilde{p_1''}$}}
    \put(12,30){\makebox(0,.2){\footnotesize$\bullet$}}
    \put(16,32){\makebox(0,3){\footnotesize$\widetilde{p_1'}$}}
    \put(21,31){\makebox(0,1){\footnotesize$\bullet$}}
    \put(21,31){\makebox(3,5){\footnotesize$\widetilde{p_0''}$}}
    \put(25,29){\makebox(0,1){\footnotesize$\bullet$}}
    \put(26,31){\makebox(0,1){\footnotesize$\widetilde{p_0'}$}}

    \put(60,32){\makebox(0,-.5){\footnotesize$\bullet$}}
    \put(55,31){\makebox(0,3){\footnotesize$\widetilde{p_1'}$}}
    \put(56,30){\makebox(0,.2){\footnotesize$\bullet$}}
    \put(60,32){\makebox(0,3){\footnotesize$\widetilde{p_1''}$}}
    \put(65,31){\makebox(0,1){\footnotesize$\bullet$}}
    \put(65,31){\makebox(3,5){\footnotesize$\widetilde{p_0'}$}}
    \put(69,29){\makebox(0,1){\footnotesize$\bullet$}}
    \put(70,31){\makebox(0,1){\footnotesize$\widetilde{p_0''}$}}

}
    \put(10,0){\footnotesize Case $T_i=L$}
    \put(60,0){\footnotesize Case $T_i=R$}
  \end{picture}

  \caption{The case of two intersection points. In the case $T_i=L$ we
  depicted two lifts of $\beta_T$, in the case $T_i=R$ we depicted
  only the segments in the grafting region.}
  \label{fig:morethanone}
\end{figure}
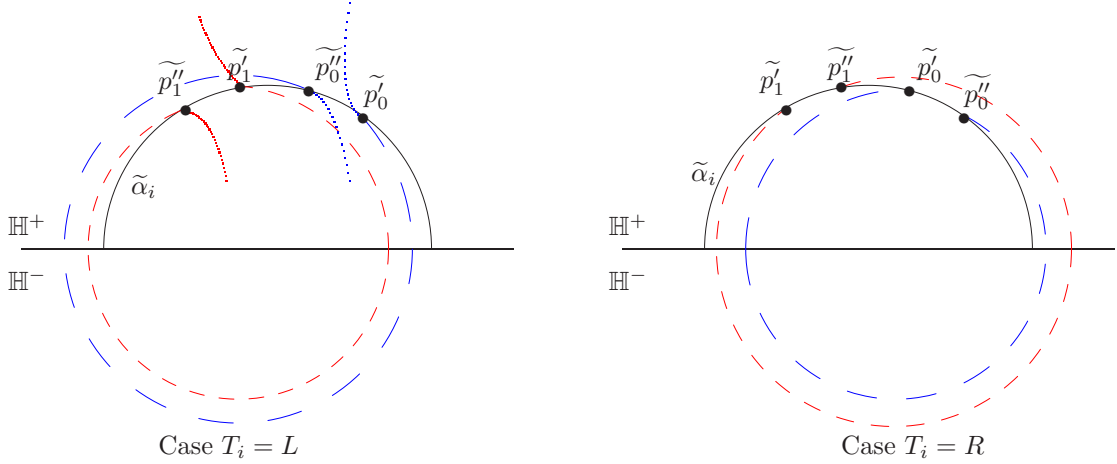

Let $\{p_j\}$ be the set of points of intersection between $\alpha_i$ and $\beta$, and  form the points $p_j'$ and $p_j''$ as before (choosing $\varepsilon$
small enough). If $\widetilde\alpha_i$ is a lift of $\alpha_i$, we see lifts $\widetilde{p_j'}$ and $\widetilde{p_j''}$ of such points. We remark that for $j\neq
k$, the point $\widetilde{p_j'}$ correspond to a lift of $\beta$ different from that of $\widetilde{p_k'}$. This is because $\alpha$ and $\beta$ are
homotopically transverse. It is worth noting at this point that it  happens that the developed images of two such lifts intersect,
but this is not a problem for our construction. Indeed, for $\beta_T$ to be graftable in $\text{Gr}_\alpha$, we only need that any single lift of $\beta_T$
is developed injectively. In Figure~\ref{fig:morethanone} we have
drawn in red (small dashed line) and blue (big dashed line)
two different lifts of $\beta_T$ entering in the same grafting region
$\mathbb{CP}^1\setminus\overline{\widetilde\alpha_i}$. The intersections of the two lifts with the grafting region are two disjoint segments, and it is clear that
such segments remain disjoint and embedded when projecting to the grafting annulus. Thus, $\beta_T$ is embedded and homotopic to $\beta$ also when multiple
intersections arise.

Let us check that any lift of $\beta_T$ develops injectively. We choose
a lift of $\beta$ and the corresponding lift of $\beta_T$. Say the red
(small dashed) lift. Since
$\alpha$ and $\beta$ are homotopically transverse, the red lift of $\beta$ intersects any lift  of any component $\alpha_i$ of $\alpha$ at most once. Thus, when
the red $\widetilde\beta_T$ enters the grafting region $\mathbb{CP}^1\setminus \overline{\widetilde\alpha_i}$, the situation is exactly that of
Figure~\ref{fig:imagebeta}. By construction, the developed image of the red $\widetilde\beta_T$ stay close to $\widetilde\alpha_i$ and its analytic prolongation
to $\mathbb H^-$. Since the lift $\widetilde\alpha_i$ is disjoint from the other lifts of $\alpha_i$ and from the lifts of different components of $\alpha$, for
$\varepsilon$ small enough the developed image of the red $\widetilde\beta_T$ is embedded.

\medskip

We now explain the variation of the construction when some $\alpha_i$ appear with multiplicity $d_i$. As was said before, it is then very important that parallel
curves have the same turning directions. In this case the grafting regions are branched coverings of $\mathbb{CP}^1$. More precisely, the universal cover of the
surface $\text{Gr}_\alpha$ is obtained by cutting $\widetilde S_u$ along the lifts $\widetilde\alpha_i$ and then by gluing back a branched covering of
$\mathbb{CP}^1$ of degree $d_i$, branched at the endpoints of $\widetilde\alpha_i$, and cut along a pre-image of $\widetilde\alpha_i$.

For any intersection point between $\alpha_i$ and $\beta$, we consider a sequence of points $p_0= p', p_1 , \ldots , p_{d_i}= p''$ in $\widetilde{\alpha_i}$
increasing from $p'$ to $p''$, and we iterate a construction similar to that of the case of multiplicity $1$. (See Figure~\ref{f3} for the situation in $S_u$ and
Figure~\ref{fig:multiple} for the situation in the universal cover.)

\setlength{\unitlength}{1ex} \begin{figure}[htbp]
  \centering
  \begin{picture}(80,20)
\multiput(0,0)(0,15){2}{\line(1,0){80}} \multiput(0,0)(20,0){5}{\qbezier(0,0)(5,7.5)(0,15)} \multiput(20,0)(20,0){4}{\bezier{20}(0,0)(-5,7.5)(0,15)}
\qbezier(0,0)(-5,7.5)(0,15)

\put(22.5,7.5){\makebox(0,0){$\bullet$}} \put(22.5,7.5){\makebox(3,-1){$p$}} \put(62.5,7.5){\makebox(0,0){$\bullet$}} \put(62.5,7.5){\makebox(-3,1){$p$}}
\put(22.3,10){\makebox(0,0){$\bullet$}} \put(22.3,10){\makebox(8,3){$p'=p_0$}} \put(62.3,5){\makebox(0,0){$\bullet$}} \put(62.3,5){\makebox(10,-3){$p''=p_2$}}
\put(42.5,7){\makebox(0,0){$\bullet$}} \put(42.5,7){\makebox(3,3){$p_1$}}

\multiput(2.5,7.5)(60,0){2}{\line(1,0){20}} \multiput(12.5,5)(60,5){2}{\makebox(0,0){$\beta$}} \qbezier(2.5,8)(15,8)(22.3,10) \qbezier(22.3,10)(30,7.5)(42.5,7)
\qbezier(62.3,5)(50,7.5)(42.5,7) \put(20,0){\qbezier(42.3,5)(45,7)(62.5,7)} \multiput(12.5,10)(60,-5){2}{\makebox(0,0){$\beta_T$}}
\put(32,9){\makebox(1,1){$\beta_T$}} \put(52,5){\makebox(0,0){$\beta_T$}}

\put(20.9,15.3){\line(3,1){6}} \put(59.1,15.3){\line(-3,1){6}} \put(40,18){\line(0,-1){2.5}} \put(10,19){\parbox{70ex}{Three copies of $\alpha_i$ bounding two
    consecutive grafting annuli }}

\put(5,12){$C_1$} \put(55,12){$C_2$}
  \end{picture}
  \caption{The curve $\beta_T$ in the surface $S_u$ when $\alpha_i$
    has multiplicity $2$. Here $T_i=L$.}
  \label{f3}
\end{figure}
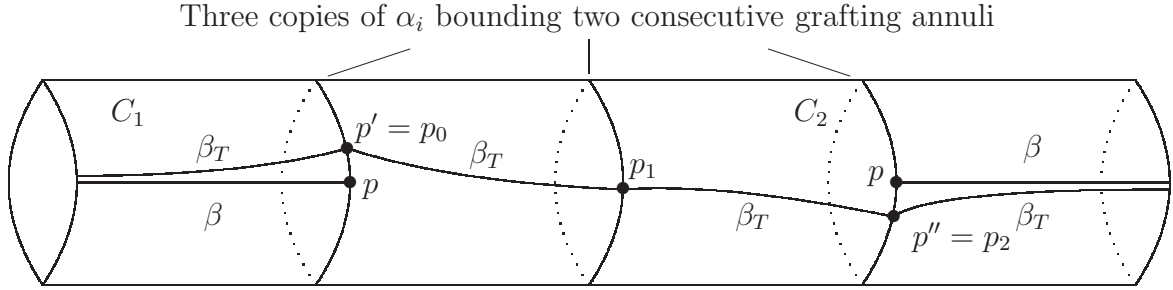

\setlength{\unitlength}{1ex} 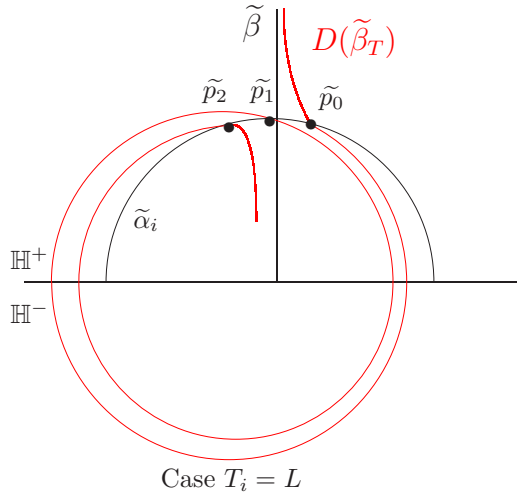
\begin{figure}[htbp]
  \centering
  \begin{picture}(30,35)
\put(0,-5){
    \put(0,20){\line(1,0){36}}
    \put(18.5,20){\line(0,1){20}}
    \put(16,38){$\widetilde\beta$}
    \put(18,20){\arc{24}{3.14}{0}}
    \put(-1,21){\footnotesize $\mathbb H^+$}
    \put(-1,17){\footnotesize $\mathbb H^-$}
{\color{red}
            \put(15,20){\arc{26}{-1.1}{3.14}}
            \put(14.5,20){\arc{25}{3.14}{0}}
            \put(15.5,20){\arc{23}{0}{4.65}}

            \qbezier(21,31.5)(19,35)(19,40)
            \qbezier(15,31.5)(17,32)(17,24.5)

            \put(21,37){$D(\widetilde\beta_T)$}
          }

    \put(8,24){\footnotesize$\widetilde\alpha_i$}

    \put(15,31.5){\makebox(0,-.5){\footnotesize$\bullet$}}
    \put(15,32.5){\makebox(-2,3){\footnotesize$\widetilde{p_2}$}}
    \put(18,32){\makebox(0,-.5){\footnotesize$\bullet$}}
    \put(18,34){\makebox(-1,0){\footnotesize$\widetilde{p_1}$}}
    \put(21,31){\makebox(0,1){\footnotesize$\bullet$}}
    \put(21,31){\makebox(3,5){\footnotesize$\widetilde{p_0}$}}

}
    \put(10,0){\footnotesize Case $T_i=L$}
  \end{picture}

  \caption{The case where $\alpha_i$ has multiplicity two. The
    grafting region is a branched covering of degree two, and $\beta_T$
  must complete two laps before exiting the region.}
  \label{fig:multiple}
\end{figure}

Finally, if some component $\beta_j$ of $\beta$ comes with multiplicity $e_j$, then we do the construction above for one copy of $\beta_j$ and then we replace
the result with $e_j$ parallel copies of the corresponding component of $\beta_T$.

\begin{rem} Note that in particular, we proved that, if $\sigma$ is a projective structure on a marked surface $S$ with Fuchsian holonomy, and $\beta$ is
\textit{any} multi-curve without component homotopic to a point, then it is possible to find a multi-curve which is graftable in $\sigma$ and isotopic to
$\beta$. It would be interesting to find conditions on a multi-curve $\beta$ that generalize  the statement for a  general projective structure (not necessarily with Fuchsian holonomy).\end{rem}

\begin{rem} There are other ways of finding graftable curves in the isotopy class of $\beta$, obtained by fixing a letter to each equivalence class of parallel
curves of the multi-curve $\beta$, instead of $\alpha$. However, this construction of multi-curve will not be discussed here. \end{rem}

\subsection{The Operation $*_T$ on homotopically transverse
  multi-curves}

Given the data $(\alpha,\beta,T)$ as in~\S\ref{ss:construction}, we produce a new isotopy class of a
multi-curve $\gamma$ in the hyperbolic  surface $S_u$ in the following
way: at each point
of intersection $p\in\alpha_i\cap\beta_j$ choose a disc $D_p$ centered
at $p$. After an isotopy we can suppose that this disc is parametrized
by an orientation
preserving map of the unit disc in the plane to $S_u$ and the image of
$\alpha_i$ corresponds to the horizontal axis and that of $\beta_j$ to
the vertical axis.

\begin{figure}
\centering
\def\svgwidth{\columnwidth}
\includegraphics[width=\textwidth]{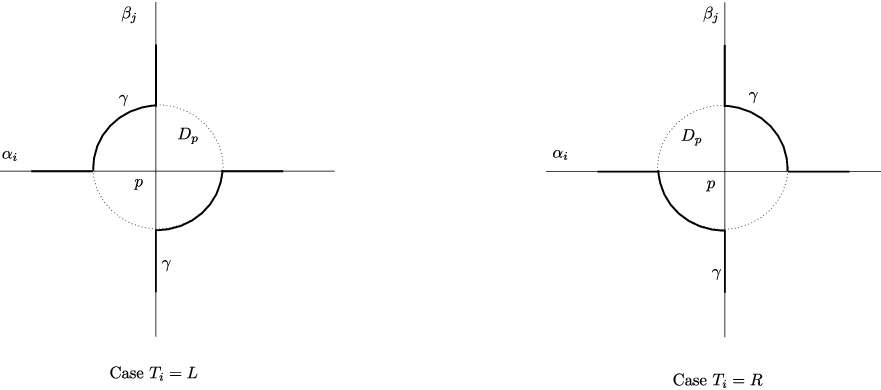}
  \caption{Construction of $\gamma$ around a point of intersection
    between $\alpha_i$ and $\beta_j$}\label{fig:intersectionalphabeta}
\end{figure}

 On $S_u\setminus\cup D_p$ the multi-curve $\gamma$ has the same
 components as $\alpha\cup\beta$. To get a multi-curve we need to join
 the endpoints by paths on $\cup\partial D_p$ by the rule given by
 $T$. As we approach an endpoint of $\alpha_i\cap \partial D_p$ from outside $D_p$ we choose the segment of $\partial D_p$ lying on the side of $\alpha_i$ given by $T_i$ between the chosen endpoint and the next point of $\beta_j\cap\partial D_p$ (see Figure~\ref{fig:intersectionalphabeta} for the two possibilites).

 This produces a family of disjoint simple closed curves $\gamma$ in
 $S_u$. The transversality condition guarantees that none of its
 components is homotopically trivial in $S_u$ and hence $\gamma$ is a
 multi-curve (see references~\cite{ITO,L}).
In the sequel, for any $(\alpha, \beta, T)$  we will denote by
$\alpha\ast_T\beta$ the resulting multi-curve:
$\alpha\ast_T\beta:=\gamma$.

\subsection{Computation of grafting annuli}
Recall that for a graftable multi-curve $\alpha$ in $S_u$ we use the
notation $\text{Gr}_\alpha=\text{Gr}(\sigma_u,\alpha)$.

\begin{prop}\label{p:construction}
Given two homotopically transverse multi-curves $\alpha$ and $\beta$,
and a set of turning directions $T$ for $\alpha$, let $\beta_{T}$
denote the graftable multi-curve constructed in~\S\ref{ss:construction}, and
$\gamma=\alpha\ast_{T}\beta$. Then
\[  \mathrm{Gr} (\mathrm{Gr}_\alpha , \beta_{T} ) =\mathrm{Gr}_\gamma.  \]
\end{prop}

\begin{proof} We have to compute the negative annuli for the structure
$\sigma' = \mathrm{Gr} (\text{Gr}_\alpha , \beta_{T} )$ given by Goldman's
theorem (see~\S\ref{ss:Goldman's parametrization}). To this
end, we will
construct a curve $\gamma_j$ in each negative annulus, and then show
that the collection of the constructed curves $\cup \gamma_j$ is
isotopic to the (graftable) multi-curve $\gamma$. By the discussion on~\S\ref{ss:Goldman's parametrization} we conclude that
$\sigma'=\text{Gr}_{\gamma}$.

First of all, note that by arguing inductively on the number of components of $\beta$, we can reduce to the case where $\beta$ is a simple loop.

To begin with, we orient $\beta$, we choose one of its lifts $\widetilde\beta$, and we number the lifts of the components of $\alpha$ that meet $\widetilde\beta$
in order of intersection with $\widetilde\beta$ as $\{\widetilde\alpha_i: i\in\mathbb Z\}$. So $\widetilde\beta$ meets $\widetilde\alpha_i$, then
$\widetilde\alpha_{i+1}$, and so on.

If $(S,\sigma)$ denotes the projective surface corresponding to the structure $\sigma=\text{Gr}_\alpha$, $\widetilde{S}$ is constructed by gluing to
$\widetilde{S_u} \setminus \bigcup \widetilde{\alpha}$ the grafting regions
 $\mathbb {CP}^1\setminus
\overline{\widetilde{\alpha}}$ (here $\widetilde\alpha$ varies among all
 lifts of all components of $\alpha$).
 Such sets will be referred to as
{\bf bubbles}. See Figure~\ref{f4}.

\setlength{\unitlength}{1ex} 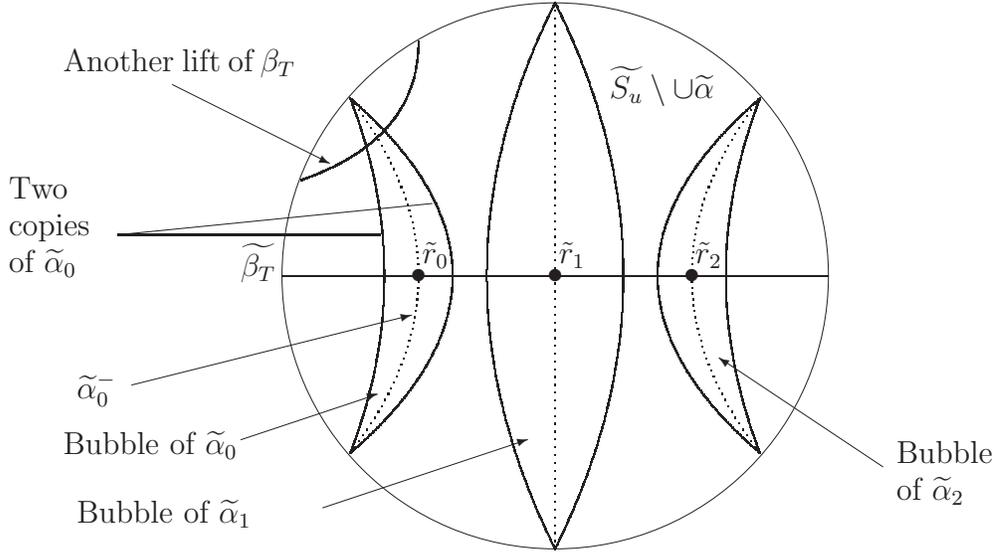
\begin{figure}[htbp]
  \centering
  \begin{picture}(80,40)
\put(40,20){\circle{40}} \qbezier(25,33)(30,20)(25,7) \qbezier(25,33)(40,20)(25,7) \bezier{60}(25,33)(35,20)(25,7)

\qbezier(40,40)(30,20)(40,0) \qbezier(40,40)(50,20)(40,0) \bezier{70}(40,40)(40,20)(40,0)

\qbezier(55,33)(50,20)(55,7) \qbezier(55,33)(40,20)(55,7) \bezier{60}(55,33)(45,20)(55,7)

\put(20,20){\line(1,0){40}} \qbezier(21.4,27)(30,30)(30,37.15)

\put(17,20){$\widetilde{\beta_T}$} \put(4,35){\parbox{17ex}{Another lift of $\beta_T$}} \put(12,34){\vector(2,-1){11.5}} \put(4,7){\parbox{13ex}{Bubble of
$\widetilde\alpha_0$}} \put(17,8){\vector(3,1){10}} \put(5,11){$\widetilde{\alpha}_0^-$} \put(9,12){\vector(4,1){20.5}} \put(0,23){\parbox{8ex}{Two copies of
$\widetilde{\alpha}_0$}} \put(8,23){\line(1,0){19.2}} \put(8,23){\line(10,1){23}}

\put(5,2){\parbox{16ex}{Bubble of $\widetilde\alpha_1$}} \put(23,3){\vector(3,1){15}} \put(65,5){\parbox{9ex}{Bubble of $\widetilde{\alpha}_2$}}
\put(64,6){\vector(-3,2){12}}

\put(44,33){$\widetilde{S_u}\setminus\cup\widetilde\alpha$}

\multiput(30,20)(10,0){3}{\makebox(0,0){$\bullet$}}
\put(30,20){\makebox(2.5,3){$\tilde r _0$}}
\put(40,20){\makebox(2.5,3){$\tilde r_1$}}
\put(50,20){\makebox(2.5,3){$\tilde r_2$}}

  \end{picture}
  \caption{The curve $\beta_T$ in $\widetilde{S}$. The bubbles
    corresponding to three consecutive lifts of components of $\alpha$ are
    depicted as ``banana'' sectors.}
  \label{f4}
\end{figure} Note that in case some component of $\alpha$ has multiplicity, then the corresponding bubbles are adjacent (this case is not depicted in the
picture).

In each bubble, let $\widetilde{\alpha}^-_i$ be the geodesic in
$\mathbb H^-$ which is the continuation of the geodesic
$\widetilde{\alpha_i}$ as a round circle
of the Riemann sphere (the dotted lines in Figure~\ref{f4}). The curve
$\widetilde{\beta_T}$ intersects these geodesics successively. For
each $n$, we denote by
$\tilde r_n$ the point of intersection of $\widetilde{\beta_T}$ and
of $\widetilde{\alpha}_n^-$. Recall that $\gamma=\alpha\ast_T\beta$
and note that  by construction $\widetilde{\gamma}$ is equivariantly
homotopic to
$\widetilde{\alpha}\ast_T\widetilde{\beta}_T$. On the other hand
$\widetilde{\alpha_i}$ is homotopic to $\widetilde{\alpha}^-_i$. A
local argument shows that
$\widetilde{\alpha}\ast_T\widetilde{\beta}_T$ is equivariantly homotopic to
$\widetilde{\alpha}^-\ast_T\widetilde{\beta}_T$. If we show that this multi-curve
is homotopic to a union of curves $\cup\gamma_j$ contained in the
negative part of $\sigma'$, and such that each connected component of
the negative part contains one of the $\gamma_j$'s we will be done.
Let us analyze the structure
$\text{Gr}(\text{Gr}_\alpha,\beta_T)$ in detail. To obtain it we have
to cut $\widetilde S$ along $\widetilde\beta_T$ and glue back a copy
of $\mathbb{CP}^1\setminus D(\widetilde\beta_T)$, where $D$ is the
developing map for $\sigma$.
Once we have cut, we have two copies $\widetilde\beta_T^R$ and
$\widetilde\beta_T^L$  of $\widetilde\beta_T$: $\widetilde\beta_T^R$
is the boundary component
that has the bubble of $\widetilde\beta_T$ on its right. In other
words, $\widetilde\beta_T^L$ is the component which is oriented
according to the orientation of
$\partial(\mathbb{CP}^1\setminus D(\widetilde\beta_T))$ .  Let $\tilde
r_n^R$ and $\tilde r_n^L$ be the points corresponding to $\tilde r_n$
lying in $\widetilde\beta_T^R$ and $\widetilde\beta_T^L$
 respectively. See these objects in Figure~\ref{f7}.

 \setlength{\unitlength}{1ex} \begin{figure}[htbp]
  \centering
  \begin{picture}(40,14)
    \qbezier(0,7)(20,18)(40,7)
    \qbezier(0,7)(20,-4)(40,7)

    \put(10,11){\makebox(0,0){$\bullet$}}
    \put(9,13){\makebox(0,0){$\tilde r_0^R$}}
    \put(10,3){\makebox(0,0){$\bullet$}}
    \put(9,1){\makebox(0,0){$\tilde r_0^L$}}

    \put(30,11){\makebox(0,0){$\bullet$}}
    \put(31,13){\makebox(0,0){$\tilde r_2^R$}}
    \put(30,3){\makebox(0,0){$\bullet$}}
    \put(31,1){\makebox(0,0){$\tilde r_2^L$}}

    \put(20,12.5){\makebox(0,0){$\bullet$}}
    \put(20,14.5){\makebox(0,0){$\tilde r_1^R$}}
    \put(20,1.5){\makebox(0,0){$\bullet$}}
    \put(20,-0.5){\makebox(0,0){$\tilde r_1^L$}}

    \qbezier(10,11)(10,11)(20,1.5)
    \qbezier(20,12.5)(24,5)(30,11)

    \qbezier(10,3)(8,6)(6,6)
    \qbezier(30,3)(33,5.5)(36,7)
\put(1,10){$\widetilde\beta_T^R$} \put(37,2.5){$\widetilde\beta_T^L$}
  \end{picture}
  \caption{The bubble of $\widetilde\beta_T$. Here the segments
    $\tilde r^R_0\tilde r^L_1$ and $\tilde r_i^R\tilde r_2^R$
    correspond to those constructed
    along the proof, contained in the negative part in the particular
    case  $T_0=L$, $T_1=L$, $T_2=R$.}
  \label{f7}
\end{figure}
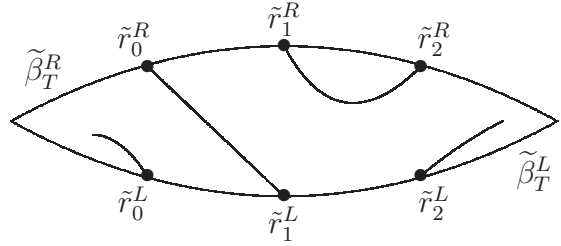

 The union of curves $\cup\gamma_j$ that we are going to
 describe in the negative part of $\sigma'$ is  a concatenation of two types of
 geodesic segments with respect to the hyperbolic metric in the
 negative part: segments contained in $\alpha^-_i$ and geodesic
 segments contained in the bubble of $\beta_T$ joining a point
$\tilde r_n^L$
 (resp. $\tilde r^R_n$) with  one of
 $\tilde r^L_{n+1},\tilde r^R_{n+1},
 \tilde r^L_{n-1},\tilde r^R_{n-1}$. The choice  will be uniquely
 defined by the sequence
 of turnings described by $T$ along $\beta_T$.  Some examples are sketched on
 Figure~\ref{f7}. These segments are most
 easily defined by using the developed image of $\widetilde{\beta}_T$
 by the developing map $D$ of $\sigma$.
 As the developed image of the points $\tilde
 r_n$  lie in the lower half plane, we can consider the geodesic
 segments joining $D(\tilde r_n)$ with $D(\tilde r_{n+1})$ for all
 $n$. Now as we cut $\mathbb{CP}^1$ along the oriented curve
 $D(\widetilde{\beta}_T)$ we realize that the pairs of points
 corresponding to each  $D(\tilde r_n)$ on each side of the cut are connected
 by the constructed segments. It is clear that for each $n$ one of the
 points in the corresponding pair is joined by a segment to one of the
 points in the pair corresponding to $D(\tilde r_{n+1})$ and the other
 to one of the points corresponding to $D(\tilde r_{n-1})$. The actual
 correspondence depends on the sequence of turnings. If $T_n=R$
 (resp. $T_n=L$) then it is $\tilde r_n^L$ (resp. $\tilde r_n^R$) that
 is joined to one of $\tilde r_{n+1}^L,\tilde r_{n+1}^R$, and this
 information is enough to determine which segments appear. Namely, if $T_n=T_{n+1}$, then the segment corresponding to
 $D(\tilde{r}_n)D(\tilde{r}_{n+1})$ describes a segment joining the
 two \emph{different} sides of the cut along $D(\tilde{\beta}_T)$. If
 $T_n\neq T_{n+1}$, the segment joins two points on the same side of
 the cut. The different possibilities before cutting
 $D(\tilde{\beta}_T)$ are sketched in Figure \ref{f:connections}.
 \begin{figure}[htbp]
  \centering
\includegraphics[width=\textwidth]{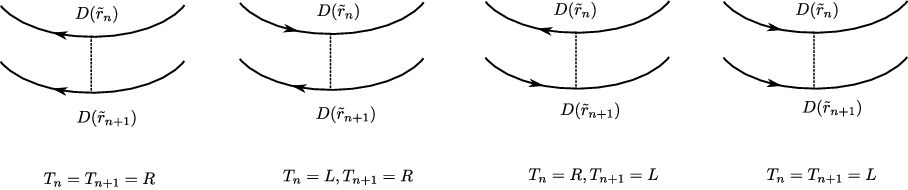}
\caption{The oriented lines represent segments of $D(\tilde{\beta}_T)$ before
  cutting. The dashed lines segments of geodesic in the negative part.}
  \label{f:connections}
\end{figure}

 After cutting $\mathbb{CP}^1$ along $D(\tilde{\beta}_T)$ we get a disc
 bounded by the two sides of the cut, that we identify with
 $\tilde{\beta}_T^L$ and $\tilde{\beta}_T^R$. Apart from that we have produced a union of disjoint segments in the disc each having one endopoint in $\{\tilde r_n^L,\tilde r_n^R\}$ and the other in $\{\tilde r_{n+1}^L,\tilde r_{n+1}^R\}$ (see
 Figure~\ref{f7} for an example of the segments obtained after the cut). The constructed segments
 produce by concatenation with those of $\tilde{\alpha}_n^-$ a union
 of curves $\cup\gamma_j$ contained in the negative part. To construct
 a homotopy with $\alpha^-\ast_T\tilde{\beta}_T$,  for each $n$ we
 choose $a_n^L$ and $a_n^R$ points on $\alpha^-_n$ lying close to
 $\tilde r_n^L$ and $\tilde r_n^R$ respectively. Remark that a segment
 in $\gamma_j$
 joining two consecutive points of the $a_n$'s has the property that
 either it cuts a single side of the cut (if the $R,L$-labels of
 $\tilde r_n$ and $\tilde r_{n+1}$ are different) or it cuts both
 sides. If it intersects only
 one side of the cut, we can homotope it with fixed endpoints to a
 segment that does not intersect the cut. Otherwise, we are obliged to
 intersect it. In fact this property characterizes the homotopy type
 with fixed endpoints of the segment. On the other hand $\gamma$ has
 the property that a segment between two consecutive $a_n$'s either
 cuts $\beta$ once (if $T_n=T_{n+1})$ or it is homotopic to a segment
 that does not intersect $\beta$ (if $T_n\neq T_{n+1})$. Therefore the
 segments between two consecutive points among the $a_n$'s of
 $\cup\gamma_j$ and $\gamma$ are homotopic with fixed endpoints. On
 the other parts of $\gamma_j$ they are equal. Therefore we can
 construct a homotopy between $\cup \gamma_j$ and
 $\widetilde{\alpha}^-\ast_T\widetilde{\beta}_T$ and the result follows.
\end{proof}

\begin{rem}
  \label{RBRBL}
Note that a corollary of Proposition~\ref{p:construction} is that
there exist graftable curves that are isotopic as curves but that
produce different structures when grafted. Indeed, let $\alpha$ and
$\beta$ two simple geodesics in the uniformizing structure such that
they intersect only in one point. Then, $\beta_R$ and $\beta_L$ are
isotopic curves (both are isotopic to $\beta$) and both graftable in
$\mathrm{Gr}_\alpha$.  By
Proposition~\ref{p:construction} we have that
$\mathrm{Gr} (\mathrm{Gr}_\alpha , \beta_{R} )
=\mathrm{Gr}_{\alpha\ast_R\beta}$ and
$\mathrm{Gr} (\mathrm{Gr}_\alpha , \beta_{L} )
=\mathrm{Gr}_{\alpha\ast_L\beta}$, which are different exotic
structures because  $\alpha\ast_R\beta$ and
$\alpha\ast_L\beta$ are not isotopic (they are positive and negative
Dehn twist of $\beta$ along $\alpha$). As the referee of this paper observed, this phenomenon was already present in Ito's work (see \cite{ITO}, Theorem 1.3).
\end{rem}

\section{Positive connectedness}

In this section we prove Theorem~\ref{t:positive connectivity}. We
begin by the following lemma, which shows that the operation $\ast_T$
is invertible.

\begin{lemma}\label{l:logarithm} Let $\alpha$ and $\gamma$ be two
  multi-curves in $S$ intersecting transversally in the sense
  of~\S\ref{ss:construction}. Suppose
that every component of $\alpha$ intersects  $\gamma$ and vice
versa. Let $T$ be a set of turning directions for $\alpha$.
 Then there exists a multi-curve $\beta$
intersecting $\alpha$ transversally in the sense
of~\S\ref{ss:construction} and such that the multi-curve
$\alpha\ast_T\beta$ is isotopic to $\gamma$.
\end{lemma}

\begin{proof} The proof is done by first constructing a multi-curve
  $\gamma'$ isotopic to the multi-curve $\gamma$ which almost
  self-intersects in a suitable way. More precisely, for each component
  $\alpha_i$ of $\alpha$, deform $\gamma$ in a small annular
  neighborhood of $\alpha_i$ as indicated in Figure
  \ref{fig:betafromgamma}, depending on the
specified turning direction. Then define the multi-curve $\beta$ as
indicated in Figure \ref{fig:betafromgamma}. It has the required
properties.
\begin{figure}
\centering
\includegraphics[width=\textwidth]{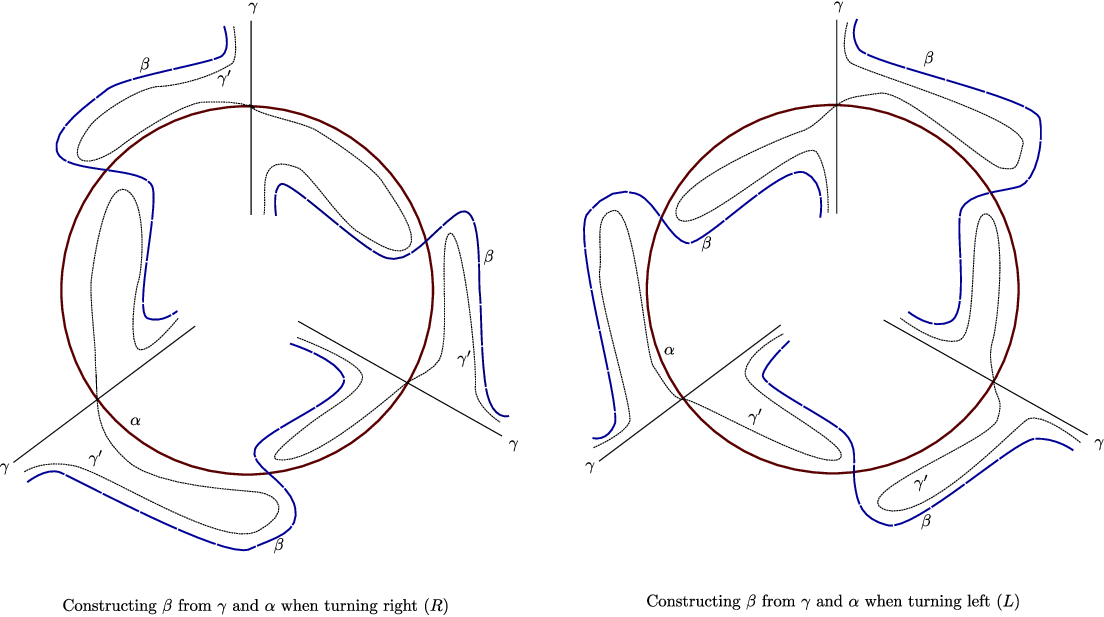}
  \caption{Constructing $\beta$}\label{fig:betafromgamma}
\end{figure}
\end{proof}

\begin{cor}
  \label{l2c}
There exists a cycle of length $2$ in the graph of multi-graftings.
\end{cor}

\begin{proof} Two symmetric applications of Lemma~\ref{l:logarithm}
  produces curves $\beta_1$ and $\beta_2$ so that
  $\mathrm{Gr}_\gamma=\mathrm{Gr}(\mathrm{Gr}_\alpha,\beta_1)$ and
$\mathrm{Gr}_\alpha=\mathrm{Gr}(\mathrm{Gr}_\gamma,\beta_2)$, proving
the existence of oriented cylces of length two.\end{proof}

We are now in a position to prove
Theorem~\ref{t:positive connectivity}. Let $(S_i,\sigma_i)$, $i=1,2$,
be projective
structures with holonomy $\rho$, both different from the uniformizing
structure $\sigma_u$. We denote by $\alpha_1$ and $\alpha_2$ the two
multi-curves coding
the negative annuli of $\sigma_1$ and $\sigma_2$ (that we think as a
multi-geodesic with multiplicities) so that
$\sigma_i=\text{Gr}_{\alpha_i}$.  Consider a
simple closed geodesic $\gamma$ cutting all components of $\alpha_1$
and all components of $\alpha_2$, and denote $(S_3,\sigma_3)=
\mathrm{Gr}_\gamma$. By two
applications of Proposition~\ref{p:construction} and
Lemma~\ref{l:logarithm}, there exist a multi-curve
$\hat{\beta_1}\subset S_1$ and a multi-curve
$\hat{\beta_2} \subset S_3$ such that $\mathrm{Gr} (\sigma_1,
\hat{\beta_1}) = \sigma_3$ and $\mathrm{Gr}(\sigma_3, \hat{\beta_2})=
\sigma_2$. This proves the
theorem.\qed

\end{document}